\documentclass[12pt]{amsart}
\usepackage{a4wide}
\usepackage{amssymb}
\usepackage{mathrsfs}
\theoremstyle{plain}
\newtheorem{theorem}{Theorem}
\newtheorem{lemma}{Lemma}[section]
\newtheorem{proposition}{Proposition}[section]

\newcommand{\NN}{\mathbf{N}}

\newcommand{\curl}{\curvearrowright}

\newcommand{\e}{\epsilon}
\newcommand{\ZZ}{\mathbf{Z}}

\newcommand{\FF}{\mathbf{F}}

\begin{document}

\title[]{Full groups and soficity\footnote{AMS
Subject Classification:20F65, 37A20}}
\author[G. Elek]{G\' abor Elek}
\thanks{Work supported in part 
by a Marie Curie grant, TAMOP 4.2.1/B-09/1/KMR-2010-003 and 
MTA Renyi "Lendulet" Groups and
Graphs Research Group.}
\begin{abstract} First, we answer a question of Pestov, by proving that the 
full group of a sofic equivalence relation is a sofic group. Then, we give a 
short proof of the theorem of Grigorchuk and Medynets that 
the topological full group of a minimal Cantor homeomorphism is LEF. 
Finally, 
we show that for certain non-amenable groups all the generalized 
lamplighter groups are sofic. 

\end{abstract} 
\maketitle
\section{Introduction}
\subsection{Sofic groups and LEF groups}
The notion of sofic groups was introduced by Weiss \cite{Weiss} 
and Gromov \cite{Gro} (in a somewhat different form) . A group 
$\Gamma$ is sofic if for any finite set $F\subset \Gamma$ and $\e>0$
there exists a finite set $A$ and a mapping $\Theta:\Gamma\to Map(A)$ such
that (\cite{ESZ1})
\begin{itemize}
\item If $f,g,fg\in F$ then 
$d_H(\Theta(fg)-\Theta(f)\Theta(g))\leq\e\,,$
where $$d_H(\alpha,\beta)=
\frac{|\{x\in A\,\mid\,\alpha(x)\neq \beta(x)\}|}{|A|}.$$
\item If $1\neq f\in F$ then
$d_H(\Theta(f),1)>1-\e\,.$
\item $\Theta(1)=1\,.$
\end{itemize}

\noindent
All amenable and residually finite groups are sofic. It is an open question 
whether non-sofic groups exist. If we add the extra requirement that
$\Theta(fg)=\Theta(f)\Theta(g)$, then we get the class of LEF-groups
(locally embeddable into finite groups). This class of groups was introduced
by Gordon and Vershik \cite{GV}. Clearly, all residually finite groups are 
LEF. However, simple, finitely presented groups are not LEF. Nevertheless, 
by a recent result of Juschenko and Monod \cite{JM} (and Theorem \ref{t2}), 
there exist simple, finitely generated LEF-groups.
\subsection{Sofic equivalence relations}
Let $X=\{0,1\}^{\NN}$ be the standard Borel space with the natural product
measure $\mu$. Let $\Phi:\FF_\infty\curl X$ be a (not necessarily free)
Borel action of the free group of countably infinite generators
$\{\gamma_1,\gamma_1^{-1},\gamma_2,\gamma_2^{-1}\dots\}$ preserving $\mu$.
Note that $\FF_\infty=\cup^\infty_{r=1}\FF_r$, where $\FF_r$ is
the free group of rank $r$. Hence, we also have probability measure preserving
(p.m.p) Borel actions $\Phi_r:\FF_r\curl X$. 
We say that $x,y\in X$ are equivalent,
$x\sim_{\Phi} y$ if there exists $w\in\FF_\infty$, such that $w(x)=y$. Note that
slightly abusing the notation we write $w(x)$ instead of $\Phi(w)(x)$.
Thus, the action $\Phi$ represents a countable measured equivalence relation
$E_\Phi$ on $X$. Similarly, each $\Phi_r$ represents a countable measured
equivalence relation $E_{\Phi_r}$ on $X$, and 
$E_\Phi=\cup_{i=1}^\infty E_{\Phi_r}$. Each equivalence relation $E_{\Phi_r}$ 
defines a graphing 
\cite{Kech} $G_r$ on $X$:
\begin{itemize}
\item $V(G_r)=X$.
\item $(x,y)\in E(G_r)$ if $\gamma_i x=y$ or $\gamma_i y =x$ for some $i$
(so, there may be loops in $G_r$).
\end{itemize}
Observe that each component of $G_r$ is a countable graph of bounded vertex 
degrees. We label each directed edge $(x,y)$ with all the generators
mapping $x$ to $y$. Thus an edge, even a loop, may have multiple labels.

\noindent
Now let us consider transitive actions of $\FF_r$ on countable sets. If 
$\alpha:\FF_r\curl Y$ is such an action then we have a bounded degree graph
structure on $Y$ with multiple labels on the edges from the set
$\{\gamma_1,\gamma_1^{-1},\dots,\gamma_r,\gamma_r^{-1}\}$. Let $T_r$ be the
set of graphs of all countable $\FF_r$-actions with a distinguished vertex 
(the root) such that all the vertices are labeled by the 
elements of $\{0,1\}^r$. Let $G\in T_r$. We define the the $k$-ball around the
root
$x$, $B_k(x)$ as the induced subgraph on vertices of $G$ in the form
of $w(x)$, where $w\in\FF_r$ is a reduced word of length at most $k$. 
That is, $B_k(x)$ is the ball centered at $x$ of radius $k$ with respect to
the shortest path metric of $G$.
The ball
$B_k(x)$ is a finite rooted graph with edge-colors from the set
$\{\gamma_1,\gamma_1^{-1},\dots,\gamma_r,\gamma_r^{-1}\}$ and vertex labels
from the set $\{0,1\}^r$. We denote the set of all possible $k$-balls 
arising from $\FF_r$-actions by $U^k_r$. We can define a compact metric 
structure on the set $T_r$ the following way. Let $d_r(G,H)=\frac{1}{2^k}$ if
$k$ is the maximal number such that the $k$-balls around the roots of $G$
resp. $H$ are isomorphic as rooted, labeled graphs.

\vskip 0.1in
\noindent
Observe that if $\Theta:\FF_\infty\curl X$ is a p.m.p action then for each
$r\geq 1$ and $x\in X$ one can associate an element $G(\Theta,x)\in T_r$. 
Namely, the orbit graph of $x$, where the vertex labels are given by the
$X$-values, restricted on the first $r$ coordinates.
Thus, we have a Borel map $\pi_\Theta:X\to T_r$\,. For $\kappa\in U^k_r$, let
$\mu_{\Theta_r}^k(\kappa)=(\pi_\theta)_\star (\mu)(L_\kappa),$ where
$L_\kappa\subset T_r$ is the set of elements $G$ such that the $k$-ball
around the root of $G$ is isomorphic to $\kappa$. In other words,
$\mu_{\Theta_r}^k(\kappa)$ is the probability that the $k$-ball around
a $\mu$-random element of $X$ is isomorphic to $\kappa$.
Now let $\alpha:\FF_r\curl Y$ be an $\FF_r$-action
on a finite set. Then for each element $y$ of $Y$, we can associate an 
element of 
$T_r$. Namely, $Y$ itself with root $y$. Hence, we can define a probability 
distribution $\mu^{k,r}_\alpha$ on $U^k_r$. Following \cite{EL} we say that the
action $\Theta:\FF_\infty\curl X$ is sofic if for all $r\geq 1$, there exists
a sequence of finite $\FF_r$-actions $\{\alpha_n\}^\infty_{n=1}$
such that for each $k\geq 1$ and $\kappa\in U^k_r$
$$\lim_{n\to\infty} \mu^{k,r}_{\alpha_n}(\kappa)=\mu^k_{\Theta_r}(\kappa)\,.$$
In \cite{EL} the authors proved that
\begin{itemize}
\item Soficity is a property of the underlying equivalence relations. That is,
if an action $\Theta_1$ is orbit equivalent to a sofic action $\Theta_2$, then
$\Theta_2$ is sofic as well.
\item Treeable equivalence relations are sofic.
\item Actions associated to Bernoulli shifts of sofic groups are sofic.
\end{itemize}
\subsection{Full groups} Let $E(X,\mu)$ be a countable, measured equivalence
relation on a Borel set $X$ with invariant measure $\mu$.
The Borel full group of $E$ is the group $[E]_B$ of all Borel bijections
$T:X\to X$ such that for any $x\in X$, $T(x)\sim_E x$. 
We call two such bijections $T_1, T_2$ equivalent if
$$\mu(\{x\in X\,\mid\, T_1(x)=T_2(x)\})=1\,.$$
The measurable full group $[E]$ is the group formed by the equivalence classes.
Obviously, $[E]=[E]_B/N$, where $N$ is the normal subgroup
of elements in $[E]_B$ fixing almost all points of $X$.

\noindent
Now, let $T:C\to C$ be a homeomorphism of the Cantor set $C$. The topological 
full group $[[T]]$ is the group of homeomorphisms $S:C\to C$ such that
$C$ can be partitioned into finitely many clopen sets $C=\cup_{i=1}^n A_i$ such
that $S_{\mid A_i}= T^{n_i}$ for some integer $n_i$.

\subsection{Results}
Answering a question of Pestov \footnote{MR2566316-MathSciNet Review}, we prove
the following theorem.
\begin{theorem}\label{t1}
The measurable full group of a sofic equivalence relation is sofic.
\end{theorem}
\noindent
Then, we give a very short proof of a result of Grigorchuk and Medynets 
\cite{GM}.
\begin{theorem}\label{t2}
The topological full group of a minimal Cantor homeomorphism is LEF.
\end{theorem}
\noindent
Let $X$ be a countably infinite set and $\Gamma$ be a countable
group acting faithfully and transitively on $X$. Then $\Gamma$ can be 
represented by automorphisms on the Abelian group $\oplus_{x\in X}\{0,1\}$.
The groups $\oplus_{x\in X}\{0,1\}\rtimes\Gamma$ are called the 
lamplighter group of the 
$\Gamma$-action. If the action is the natural translation action 
on $\Gamma$, then
we get the classical lamplighter group of $\Gamma$.
Paunescu \cite{Pau} proved that if $\Gamma$ is sofic, then the
classical lamplighter group 
$\oplus_{\gamma\in\Gamma}\{0,1\}\rtimes \Gamma$ is sofic.
If $\Gamma$ is amenable, then all its generalized lamplighter groups
are amenable hence sofic. Nevertheless, we show that there exist
non-amenable groups for which all the generalized lamplighter groups are
sofic.
\begin{theorem}\label{t3}
Let $\Gamma^k$ be the $k$-fold free product of the cyclic group
of two elements. Then, for any transitive, faithful action of $\Gamma^k$
on a countable set the associated lamplighter group is LEF.
\end{theorem}

\noindent
{\bf Acknowledgement:} The author thanks Nicolas Monod and G\'abor Pete for 
valuable discussions. 
\section{Compressed sofic representations}
Let $\Gamma$ be a countable sofic group with elements 
$\{\gamma_1,\gamma_2,\dots\}$. A compressed sofic representation of $\Gamma$
is defined the following way. For any $i\geq 1$, we have a constant
$\e_i>0$ and for any $n\geq 1$ we have mappings $\Theta_n:\Gamma\to Map(A_n)$
such that $|A_n|<\infty$ satisfying the following condition: For all $r>0$
and $\e>0$ there exists $K_{r,\e}>0$ such that if $n>K_{r,\e}$ then
\begin{itemize}
\item $d_H(\Theta_n(\gamma_i\gamma_j)\Theta_n(\gamma_i)\Theta_n(\gamma_j))<\e$ 
if $1\leq i,j \leq r$.
\item $d_H(\Theta_n(\gamma_i),Id)>\e_i$ if $1\leq i \leq r$.
\end{itemize}
Thus, in a compressed sofic representation we allow large amount of
fixed points for each $\gamma\in\Gamma$. 
\begin{lemma}
If $\Gamma$ has a compressed sofic representation then $\Gamma$ is sofic.
\end{lemma}
\proof
Let $\tilde{\Theta}^k_n:\Gamma\to Map(A^k_n)$ be defined by
$$\tilde{\Theta}^k_n(\gamma)(x_1,x_2,\dots,x_k)=(\Theta_n(\gamma)(x_1),
\Theta_n(\gamma)(x_2),\dots)\,.$$
Observe that if $\gamma,\delta\in\Gamma$, then
\begin{itemize}
\item $d_H(\tilde{\Theta}^k_n(\gamma\delta), \tilde{\Theta}^k_n(\gamma)
\tilde{\Theta}^k_n(\delta))\leq (1-d_H(\Theta_n(\gamma\delta),
\Theta_n(\gamma)\Theta_n(\delta))^k$
\item $d_H(\tilde{\Theta}^k_n(\gamma),Id)>1-(1-d_H(\Theta_n(\gamma),Id))^k$
\end{itemize}
Hence, we can choose $\epsilon$, $n$ and $k$ appropriately to obtain for
any $F\subset\Gamma$ and $\e'>0$ a map $\Theta$ as in the Introduction, 
proving the soficity of $\Gamma$. \qed
\section{The proof of Theorem \ref{t1}}
Let $\Phi:\FF_\infty\curl\{0,1\}^{\NN}$ be a sofic action preserving
the product measure $\mu$. Let $\Gamma\subset [E]$ be a finitely
generated group, where $[E]$ is the equivalence relation defined by $\Phi$.
So, we have an action $\Phi_{\Gamma}:\Gamma\curl\{0,1\}^{\NN}$.
Our goal is to construct a compressed sofic representation of $\Gamma$.
Let $\{\gamma_n\}^\infty_{n=1}$ be an enumeration
of the elements of $\Gamma$. Let $\e_n=\mu(Fix(\Phi_\Gamma(\gamma_n))/2$.
Since $\Gamma$ is in the full group, $\e_n>0$.
Now, fix a subset $F\subseteq \Gamma$ and $\e>0$. We need to construct
a map $\Theta:F\to Map(A)$ for some finite set $A$ such that
if $\gamma_i,\gamma_j,\gamma_{i}\gamma_{j}\in F$ then
\begin{equation} \label{e1}
d_H(\Theta(\gamma_i\gamma_j)\Theta(\gamma_i)\Theta(\gamma_j))< \e
\end{equation}
\begin{equation} \label{e2}
d_H(\Theta(\gamma_i),1)>\e_i
\end{equation}

\vskip 0.2in

Let $\{s_1,s_1^{-1},s_2,s_2^{-1},\dots,s_m,s_m^{-1}\}$ be a symmetric
generating set for $\Gamma$. Observe that we have an action 
$\Sigma_\Gamma:\FF_m\curl \{0,1\}^{\NN}$ preserving $\mu$
such that $\Sigma_\Gamma(\delta)=\Phi_\Gamma(\tau(\delta))$, where
$\tau:\FF_m\to\Gamma$ is the natural quotient map. 
A dyadic $E$-map of depth $k$ is a Borel map
$Q:X\to X$ is defined the following way. For each $\rho\in \{0,1\}^k$
we pick $w_Q(\rho)\in \FF_k\subset\FF_\infty$
and define $Q(x)=\Phi(w_Q(\rho))(x)$ if the first $k$-coordinate 
of $x$ is $\rho$.

\noindent
A dyadic approximation of $\Gamma$ is a sequence of families
$\{Q_k(s_i)\}^m_{i=1}, \{Q_k(s_i^{-1})\}^m_{i=1}$, where
for any $1\leq i \leq m$
\begin{itemize}
\item $Q_k(s_i):X\to X$, $Q_n(s_i^{-1}):X\to X$ are dyadic $E$-maps
of depth $k$.
\item $\lim_{k\to\infty} \mu(\{x\in X\,\mid\, 
Q_k(s_i)(x)\neq\Sigma_\Gamma(s_i)(x)\})=0$
\item $\lim_{k\to\infty} \mu(\{x\in X\,\mid\, 
Q_k(s_i^{-1})(x)\neq\Sigma_\Gamma(s_i)(x)\})=0$
\end{itemize}

\noindent
We do not require $Q_k$ to be a bijection. Nevertheless,
$Q_k$ can be extended to a homomorphism
from $\FF_m$ to $Map(X)$. Note that for simplicity we 
identified the generating set of $\FF_m$ by the set $\{s_1,s_1^{-1}, s_2, 
s_2^{-1},\dots, s_m,s_m^{-1}\}$.

Since all the
$\Sigma_\Gamma(s_i)'s$ are Borel bijections such dyadic approximations clearly
exist. The following lemma is an immediate consequence of the definition
of the dyadic approximation.
\begin{lemma}
For any $\delta\in\FF_m$
$$\lim_{k\to\infty} \mu(Fix(Q_k(\delta)))=\mu(Fix(\Sigma_\Gamma(\delta)))\,.$$
\end{lemma}
\begin{proposition}\label{mapro}
There exists a sequence of mappings
$\hat{\Theta}_k:\FF_m\to Map(B_k)$, where $|B_k|<\infty$ such
that for any $\delta\in \FF_m$
$$\lim_{k\to\infty}(\mu(Fix(Q_k(\delta)))-
\frac{|Fix (\hat{\Theta}_k(\delta))|}{|B_k|})=0\,.$$
\noindent
That is
$$\lim_{k\to\infty}\frac{|Fix (\hat{\Theta}_k(\delta))|}{|B_k|}=
\mu(Fix(\Sigma_\Gamma(\delta)))\,.$$
\end{proposition}
\proof
Let $\Phi_k:\FF_k\curl\{0,1\}^{\NN}$ be the restriction of $\Phi$.
Since $\Phi$ is sofic, there exists a sequence of mappings
$\{\iota^n_k:\FF_k\curl Perm(C_{k,n})\}^\infty_{n=1}$, where $C_{k,n}$ is a finite
$\{0,1\}^k$-vertex labeled graph such that for any $t\geq 1$ and
$\kappa\in U^t_k$
$$\lim_{n\to\infty} \mu^{t,k}_{\iota^n_k}(\kappa)=\mu^t_{\Phi_k}(\kappa)\,.$$

\noindent
Recall that $Q_k$ is not necessarily an action, only a homomorphism from
$\FF_m$ to $Map(X)$. Hence, the local statistics of $Q_k$ can not be
described using the elements of $U^t_k$ as 
in the case of honest $\FF_m$-actions. So, let $W^t_k$ be the set of isomorphism
classes of rooted $t$-balls of vertex degrees at most $2m$, where
the vertices are labeled by elements of the set $\{0,1\}^k$ and
the edges (possibly loops) are labeled by subsets of
$\{s_1,s_1^{-1},s_2,s_2^{-1},\dots,s_m,s_m^{-1}\}$. Note that $U^t_k\subset W^t_k$.
Let $x,y\in X$ be points such that
$B^{\Phi_k}_{k^2}(x)$ and $B^{\Phi_k}_{k^2}(y)$ represent the same element in
$U^{k^2}_k$. Here $B^{\Phi_k}_{k^2}(x)$ denotes the $k$-ball with respect to
the graphing associated to $\Phi_k$. Then, by the definition of the
dyadic approximations
$B^{Q_k}_k(x)$ and $B^{Q_k}_k(y)$ represent the same elements in $W^k_k$.
Now we construct a sequence of maps $\hat{\Theta}^n_k:\FF_m\curl Map(C_{k,n})$ 
the following way.
$$\hat{\Theta}^n_k(s_i)(x)=\iota^n_k(w_{Q_k(s_i)}(\rho(x)))(x)\,,$$
where $\rho(x)$ is the $\{0,1\}^k$-label of $x$.
By the previous observation, for any $\delta\in \FF_m$
$$\lim_{n\to\infty} \frac{|Fix(\hat{\Theta}^n_k(\delta))|}{|C_{k,n}|}=
\mu(Fix(Q_k(\delta)))\,.$$
\noindent
This finishes the proof of the proposition\qed

\vskip 0.2in
\noindent
Pick a section $\sigma:\Gamma\to\FF_m$, that is a map such that
$\tau\sigma=Id$. Let $\hat{\Theta}_k$ as in Proposition \ref{mapro}.
Define $\Theta_k:\Gamma\to Map(B_k)$ by
$$\Theta_k(\gamma)=\hat{\Theta}_k(\sigma(\gamma))\,.$$
Then $\{\Theta_k\}^\infty_{k=1}$ is a compressed sofic representation
of $\Gamma$. \qed

\section{The proof of Theorem \ref{t2}}
Let $T:C\to C$ be a minimal homeomorphism and $\Gamma\subset [[T]]$ be
a finitely generated subgroup of the topological full group of $T$ 
with symmetric
generating set $S=\{a_1,a_2,\dots,a_k\}$. It is enough to prove that $\Gamma$ is
LEF. Let $x\in C$ and consider the $T$-orbit $\{T^n(x)\}^\infty_{-\infty}.$
We define the map $\phi:\Gamma\to Perm(\ZZ)$ of $\Gamma$ into the 
permutation group of the integers the following way.
Let $\phi(\gamma)(n)=m$, if $\gamma(T^n(x))=T^m(x)\,.$
Since $T$ acts freely on $C$, $\phi$ is well-defined.
\begin{lemma} \label{l14}
$\phi$ is an injective homomorphism.
\end{lemma}
\proof If $\phi(\gamma)=Id$, then $\gamma$ fixes all the elements
of the orbit of $x$. Since all the orbits are dense, this implies
that $\gamma=1$. The fact that $\phi$ is a homomorphism follows immediately,
since $\phi$ is the restriction of the $\Gamma$-action onto the orbit of $x$.
\qed

\vskip 0.2in
\noindent
Let $a=\max |n|$, where for some $p\in C$ and $a_i\in S$, $a_i(p)=T^n(p)\,.$
We define a sequence 
$$l:\ZZ\to\{-a,-a+1,\dots,0,1,\dots, a-1,a\}^S$$
the following way. 
Let $l(n):=(t_{a_1},t_{a_2},\dots,t_{a_k})$, where
$a_i(T^n(x))=T^{n+t_{a_i}}(x)\,.$
The following lemma is well-known, we prove it for the sake of completeness.
\begin{lemma}\label{l15}
$l$ is a repetitive sequence, that is, if we find a substring $\sigma$ in $l$,
then there exists $m\geq 1$ such that for any interval of length $m$ we can find
$\sigma$.
\end{lemma}
\proof
For a point $p\in C$, we can define its $n$-pattern
$$q_n(p):=\{-n,-n+1,\dots,0,1,\dots,n-1,n\}\to \{-a,-a+1,\dots,a-1,a\}$$
by $q_n(p)(j):=(t_{a_1},t_{a_2},\dots,t_{a_k})$, where
$a_i(T^j(x))=T^{j+t_{a_i}}(x)\,.$
Observe that the set of points with a given $n$ pattern is closed.
Now, let us suppose that for a sequence $\{k_r\}^\infty_{r=1}\subset\ZZ$
the intervals $(k_r-r,k_r+r)$ do not contain $\sigma$ as a substring.
Then, if $z$ is a limit point of $\{T^{k_r}(x)\}^\infty_{r=1}$, no translates
of $z$ have $\sigma$ as a part of their $n$-patterns. Therefore the orbit
closure of $z$ does not contain $x$, in contradiction with the
minimality of $T$.\qed

\vskip 0.2in
\noindent
Now let $r\geq 1$ and consider the string 
$\sigma_r=l_{\mid\{-ar,-ar+1,\dots,ar-1,ar\}}$, where $a$ is the constant
defined above. Note that if $\gamma\in \Gamma$ is the product of at most $r$
generators then $|\phi(\gamma)(i)-i|\leq ar$. Pick $n>10a^{r}$ such that
\begin{itemize}
\item $l_{\mid\{-ar+n,-ar+1+n,\dots,ar-1+n,ar+n\}}=\sigma_r$,
\item for any $\gamma\in \Gamma$ that is the product of at most $r$ generators
there is $0<j<n$ such that $\gamma(j)\neq j$.
\end{itemize}
Now we define $\phi_r:W^r\to Perm(\ZZ_n)$, where $W^r$ is the
set of elements in $\Gamma$ that are products of at most $r$ generators
by $\phi_r(i)=\phi(i)(mod\, n)$.
Clearly, $\phi_r$ is injective and if $x,y,xy\in W^r$ then 
$\phi_r(x)\phi_r(y)=\phi_r(xy)$. This implies that $\Gamma$ is LEF. \qed

\section {The proof of Theorem \ref{t3}}
Let $\alpha:\Gamma^k\to X$ be a transitive and faithful action of
the free product group. Consider the Schreier graph $G_\alpha$ of the action
with respect to the generators of the $k$ cyclic groups $\{a_1,a_2,\dots a_k\}$.
Recall that $V(G_\alpha)$ is $X$ and $(x,y)\in E(G)$ if
$y=a_ix$ for some $i\geq 1$. Hence $G_\alpha$ is a connected graph of vertex 
degree
bound $k$. 
\begin{proposition}\label{p8} Let $\alpha$ be as above. Then
for any $1\neq w\in\Gamma^k$, there exist infinitely many $y\in X$ such that
$\alpha(w)(y)\neq y$.
\end{proposition}
\proof We will need the following lemma.
\begin{lemma}
For any finite set $S\subseteq X$, there exists $g\in \Gamma^k$ such that
$gS\cap S=\emptyset$.
\end{lemma} 
\proof We define a lazy random walk on $X$ the following way. For $y\in X$
the transition probability $p(x,y)=l/k$, where $l$ is the number of generators
$a_i$ such that $a_ix=y$. It is well-known (see e.g. \cite{MP},\cite{Russ})
that the probabilities $p_n(x,y)$ tend to zero for each pair $x,y\in X$.
Now consider the standard random walk on the Cayley graph of $\Gamma^k$, the 
$k$-regular tree. Let $P_n(g)$ be the probability being at $g$ after taking
$n$ steps starting from the identity. Then,
$$p_n(x,y)=\sum_{g\in\Gamma, gx=y} P_n(g)\,.$$ By the previous observation, if
$n$ is large enough, then
$$\sum P_n(g)<1\,,$$
where the summation is taken for all $g\in\Gamma^k$ such that $gx\in S$, for 
some $x\in S$. Hence, there exists $g\in\Gamma^k$ such that
$gS\cap S=\emptyset$.\quad\qed

\vskip 0.2in
\noindent
Now let us suppose that
$w\in\Gamma^k$ fixes all points of $X$ outside a finite set $S$. That is 
$\alpha(w)(S)=S$. Let $gS\cap S=\emptyset$. Then $gwg^{-1}$ fixes all the
points of $X$ outside $gS$. Therefore the commutator $[w,gwg^{-1}]$
fixes all elements of $X$, in contradiction with the assumption that
the action is faithful. \qed

\vskip 0.1in
Now fix a vertex $x\in X$ and consider the ball of radius $n$, $B_n(x)$ around
$x$.
We define an action $\alpha_n:\Gamma^k\curl B_n(x)$ the following way.
Let $\partial B_n(x)$ be the boundary of the ball $B_n(x)$, that is, the
set of all $y\in B_n(x)$ such that there exists $a_i$ for which
$\alpha(a_i)y\notin B_n(x)$.
If $y\notin \partial B_n(x)$, then let $\alpha_n(a_i)y=\alpha(a_i)y\,.$
If $y\in \partial B_n(x)$ and $\alpha(a_i)y\notin B_n(x)$, then
let $\alpha_n(a_i)(y)=y\,.$ Finally, if 
$y\in \partial B_n(x)$ and $\alpha(a_i)y\in B_n(x)$, then
let $\alpha_n(a_i)(y)=\alpha (a_i)(y).$
Now let $L^n_k=\{0,1\}^{B_n(x)}\rtimes_{\alpha_n}\alpha_n(\Gamma^k)$ be the 
associated
finite lamplighter group and $L^k=\oplus_{x\in X}\{0,1\}\rtimes_\alpha \Gamma^k$.
Our goal is to embed $L^k$ into $L^n_k$ locally. That is, for any finite
set $F\subset L^k$ we construct an injective map $\Theta:F\to L^n_k$ such
that $\Theta(fg)=\Theta(f)\Theta(g)\,.$ Recall, that each element of
$L^k$ can be uniquely written in the form $a\cdot w$, where
$a\in \oplus_{x\in X}\{0,1\}$ and $w\in\Gamma^k$.
We regard the elements of the lamplighter group as permutations
of the set $\oplus_{x\in X}\{0,1\}$. If $\kappa\in \oplus_{x\in X}\{0,1\}$ and
$p\in X$
then
$$(a\cdot w)(\kappa)_{\mid p}=a(p)+\kappa(\alpha(w^{-1})(p))\,.$$
We will also use the product formula
$$(a_2\cdot w_2)(a_1\cdot w_1)=(a_2+\alpha(w_2)(a_1), w_2 w_1)\,,$$
where $\alpha(w_2)(a_1)(q)=a_1(\alpha(w_2^{-1})(q)).$
For $l\geq 1$, let $H_l$ be the set of elements of $L^k$ in the 
form of $a\cdot w$, where $w$ is a word of length at most $l$ and
the support of $a$ is contained in $B_l(x)$. For $n\geq l$ we define
the map $\tau^n_l:H_l\to L^n_k$ by
$\tau^n_l(a\cdot w):=a\cdot\alpha_n(w)$. 
\begin{lemma} \label{l11}
If $n$ is large enough then $\tau^n_l$ is injective.
\end{lemma}
\proof
If $n$ is large enough then $B_n(x)$ contains a point $y$ such that
\begin{itemize}
\item $\alpha(w)(y)\neq y$
\item $d(y,\partial B_n(x))>l$
\item $d(y,B_l(x))>l$,
\end{itemize}
where $d$ is the shortest path distance on the Schreier graph $G_\alpha$.
Let $\kappa\in \oplus_{x\in X}\{0,1\}$ be the element which is
$1$ at $y$ and zero otherwise. Then
$$\tau^n_l(a\cdot w)(\kappa)_{\mid \alpha_n(w)(y)}=1\,,$$ hence
$\tau^n_l(a\cdot w)$ is not trivial. \qed

\noindent
The following lemma finishes the proof of Theorem \ref{t3}.

\begin{lemma} \label{l12}
Suppose that $(a_1\cdot w_1),(a_2\cdot w_2)$ and
$(a_2\cdot w_2)(a_1\cdot w_1)\in H_l$ and $n$ is large enough. Then
$$\tau^n_l((a_2\cdot w_2))\tau^n_l((a_1\cdot w_1))=
\tau^n_l((a_2\cdot w_2)(a_1\cdot w_1))\,.$$ \end{lemma}
\proof
We need to prove that
$$(a_2\cdot \alpha_n(w_2))(a_1\cdot\alpha_n(w_1))=(a_2+\alpha(w_2)(a_1))
\cdot\alpha_n(w_2w_1)$$
holds in $L^n_k$. 
Fix an element $\kappa\in\{0,1\}^{B_n(x)}\,.$ Let $n>10l$ and
$d(p,\partial B_n(x)))>5l\,.$
Then 
$$(a_2\cdot \alpha_n(w_2))(a_1\cdot\alpha_n(w_1))(\kappa)_{\mid p}=
(a_2\cdot w_2)(a_1\cdot w_1)(\overline{\kappa})_{\mid p}$$
and
$$(a_2+\alpha(w_2)(a_1)\cdot \alpha_n(w_2w_1))(\kappa)_{\mid p}=
(a_2+\alpha(w_2)(a_1)\cdot (w_2w_1)(\overline{\kappa})_{\mid p}\,,$$
where $\overline{\kappa}$ is an extension of $\kappa$ onto $X$.
On the other hand, if
$d(p,\partial B_n(x)))\leq5l\,$, then
$$(a_2\cdot \alpha_n(w_2))(a_1\cdot\alpha_n(w_1))(\kappa)_{\mid p}=
\alpha_n(w_2)\alpha_n(w_1)(\kappa)_{\mid p}=$$
$$=\alpha_n(w_2w_1)(\kappa)_{\mid p}=
(a_2+\alpha(w_2)(a_1))\cdot \alpha_n(w_2w_1))(\kappa)_{\mid p}\,\quad\qed$$

\vskip 0.3in
\noindent
gabor.elek@epfl.ch

\end{document}